\documentclass[11pt]{amsart}
\usepackage{amsfonts,amssymb,amscd}

\newcommand{\bir}{\mathop{\rm {- - }\to}\nolimits}
\newcommand{\pal}{\text{---}}
\newcommand{\mt}[1]{\operatorname{#1}}
\newcommand{\wt}{\mt{wt}}
\newcommand{\ord}{\mt{ord}}
\newcommand{\ex}{{\mt{ex}}_4}
\newcommand{\gr}{{\mt{gr}}^0_C\omega}
\newcommand{\gro}{{\mt{gr}}^1_C\OOO}

\newcommand{\bb}[1]{{\mathbb{#1}}}
\newcommand{\cc}[1]{{\mathcal{#1}}}
\newcommand{\CC}{{\bb{C}}}

\newcommand{\ZZ}{{\bb{Z}}}
\newcommand{\QQ}{{\bb{Q}}}
\newcommand{\PP}{{\bb{P}}}
\newcommand{\NN}{{\bb{N}}}
\newcommand{\OOO}{{\cc{O}}}

\newcommand{\m}{\bar{m}}
\newcommand{\qq}{\begin{flushright}
$\square$\end{flushright}}
\newcommand{\cyc}[1]{\bb{Z}_{#1}}
\newcommand{\3}{^{\sharp}}
\newcommand{\6}{^{\flat}}

\newtheorem{theorem}[subsection]{Theorem}
\newtheorem{theorem1}[subsubsection]{Theorem}
\newtheorem{proposition}[subsection]{Proposition}
\newtheorem{proposition1}[subsubsection]{Proposition}
\newtheorem{conjecture}[subsection]{Conjecture}
\newtheorem{lemma}[subsection]{Lemma}
\newtheorem{lemma1}[subsubsection]{Lemma}
\newtheorem{corollary}[subsection]{Corollary}
\newtheorem{corollary1}[subsubsection]{Corollary}

\theoremstyle{definition}
\newtheorem{definition}[subsection]{Definition}
\newtheorem{definition1}[subsubsection]{Definition}
\newtheorem{example}[subsection]{Example}

\newtheorem{proposition-definition}[subsection]{Proposition-Definition}
\newtheorem{proposition-definition1}[subsubsection]{Proposition-Definition}

\theoremstyle{remark}
\newtheorem{remark}[subsection]{Remark}

\newtheorem*{remark3}{Remark}

\begin{document}

\title[On Extremal Contractions]
{On Extremal Contractions from Threefolds to Surfaces:\\
the Case of One non-Gorenstein Point}

\author{Yuri G. Prokhorov}

\email{
prokhoro@mech.math.msu.su\\
prokhoro@nw.math.msu.su
}
\thanks{The author was supported in part by
the Russian Foundation of Fundamental Research
\# 96-01-00820 and
the International Science Foundation \# M9O 300.}

\address{Algebra Section,
Department of Mathematics,
Moscow State University,
Vorob'evy Gory, Moscow
117 234, Russia}

\date{}

\dedicatory{In Memoriam: Professor Wei-Liang Chow}

\begin{abstract}
Let $f\colon X\to S$ be an extremal contraction from
a threefold with only terminal singularities to a surface.
We study local analytic structure such contractions near
degenerate fiber $C$  in the case when $C$  is irreducible and
$X$ has on $C$ only one non-Gorenstein point.
\end{abstract}

\maketitle
\section*{Introduction}
Let $V$ be a non-singular algebraic projective threefold over $\CC$
 of Kodaira dimension
$\kappa(V)=-\infty$. It follows from the Minimal Model Program
that $V$ is birationally equivalent to an algebraic projective
threefold $X$
with only terminal $\QQ$-factorial singularities  and with
 a fiber space structure
$f\colon X\to S$ over a lower  dimensional variety $S$ such that $-K_X$
is $f$-ample
and relative Picard number of $X/S$ is one.
Thus, at least in principle, one can reduce the problem of studying
birational
structure of $V$ to studying structure of $S$ and the fiber space $X/S$.
In the present paper we investigate the local structure of
$X/S$ near singular fibers in the case $\dim(S)=2$.
In this situation it is easy to see that  a general fiber is a conic
in $\PP^2$.
We will work in the analytic situation.
\begin{definition}\label{def}
Let $(X,C)$ be a germ of a three-dimensional complex
space   along a  compact reduced curve $C$ and let
$(S,0)$ be a germ of a
two-dimensional normal complex space.
Suppose that $X$ has  at worst terminal singularities.
Then we say that a proper morphism
$f\colon (X,C)\to (S,0)$ is a {\it Mori conic bundle}
if \par
\begin{enumerate}
\renewcommand\labelenumi{(\roman{enumi})}
\item $(f^{-1}(0))_{\mt{red}}=C$;
\item $f_*\cc{O}_X=\cc{O}_S$;
\item $-K_X$ is ample.
\end{enumerate}
 \end{definition}
The following conjecture is interesting for applications of
Sarkisov program to  the
  rationality problem for conic bundles \cite{Iskovskikh1},
  \cite{Iskovskikh}
 or studying of $\QQ$-Fano
threefolds with extremal contractions onto a surface \cite{Pro1}.

\begin{conjecture}\label{11}
Let $f\colon (X,C)\to (S,0)$ be a Mori conic bundle.
Then $(S,0)$ is a DuVal singularity of type $A_n$.
\end{conjecture}
It was proved in \cite{Cut} that if $X$ is Gorenstein, then
$(S,0)$ is non-singular and $f\colon (X,C)\to (S,0)$ is a
"usual" conic bundle,
i.~e. there exists an embedding
$(X,C)\subset\bb{P}^2\times (S,0)$ such that
 all the fibers of $f$ are conics in $\PP^2$.
Earlier \cite{Pro2} (see also \cite{Pro3}) we classified all
Mori conic bundles containing
only points of indices $\le 2$. In this case $(S,0)$ is
non-singular or
DuVal of type $A_1$.
In general,
Conjecture \ref{11} follows from the following special
case of Reid's general elephant conjecture
(see \cite{Pro}, \cite{Pro1}).
\begin{conjecture}\label{reidd}
Let $f\colon (X,C)\to (S,0)$ be a Mori conic bundle.
Then the linear system $|-K_X|$ contains a divisor having
 only DuVal singularities.
\end{conjecture}
 It can be generalized as follows.
\begin{conjecture}
\label{log}
Let $f\colon (X,C)\to (S,0)$ be a Mori conic bundle. Then for
any $n\in\NN$
the following holds
\par\medskip\noindent
$(*_n)$\qquad
There exists a divisor $D\in |-nK_X|$ such that
$K_X+\frac{1}{n}D$ is log-terminal.
\par\medskip\noindent
\textup{(}Note that \textup{\ref{log}} $(*_1)$ is equivalent
to \textup{\ref{reidd}}\textup{)}.
\end{conjecture}
 By  \cite[Theorem 4.5]{Kawamata},
for every threefold $X$ with terminal singularities
there exists a projective bimeromorphic morphism
$q\colon X^{q}\to X$ called {\it $\QQ$-factorialization} of $X$
 such that $X^{q}$ has only terminal
$\QQ$-factorial along $q^{-1}(C)$ singularities and $q$ is an
isomorphism in
codimension 1. If $f\colon (X,C)\to (S,0)$ is a Mori conic bundle,
then applying the Minimal Model Program to $X^{q}$ over $(S,0)$
we obtain a  new Mori conic bundle $f'\colon (X^{\prime},C')\to (S,0)$
over the  same base surface $S$ and
$\rho (X^{\prime},C^{\prime})/(S,0)$=1.
In particular  $C'$ is irreducible by
Corollary~\ref{vanishing}.
Therefore it is sufficient to prove Conjecture~\ref{11}
assuming that $C$ is irreducible.
\par
In the this paper we will study Mori conic bundles
$f\colon (X,C)\to (S,0)$ with irreducible
central fiber. In this situation $X$ contains at most three
singular points (Corollary~\ref{<=3}). Arguments similar to
\cite[0.4.13.3]{Mori-flip} show also that $X$ can contain  at
most two
 non-Gorenstein points
\cite{Pro1}.  Our results concern only the case when $X$ contains
  only one non-Gorenstein singular point.
It will be proved in this case (Theorem~\ref{main}) that
 the point $(S,0)$ may be either
 non-singular or DuVal of type $A_1$ or $A_3$. We also show  that
either  $(*_1)$ or $(*_2)$ hold and give a
classification if $(S,0)$ is singular.
Our main  tool will be Mori's technique \cite{Mori-flip}.
\subsection*{Acknowledgments.}
I would like to acknowledge discussions that I had on this
subject with Professors V.I.~Iskovskikh,
V.V.~Shokurov and M.~Reid.
I have been working on this problem at
Max-Planck Institut f\"ur Mathematik and
the Johns Hopkins University. I am very grateful them
for hospitality and support.

\section{Preliminary results}

The following is a consequence of the Kawamata-Viehweg
vanishing theorem (see \cite[\S 4]{Nakayama},
 \cite[1-2-5]{KMM}).

\begin{theorem}
Let $f\colon (X,C)\to (S,0)$ be a
Mori conic bundle. Then
$R^if_*\cc{O}_X=0$, $i>0$.
\qq
\end{theorem}

\begin{corollary1}
\textup{(}cf. \cite[(1.2)-(1.3)]{Mori-flip}\textup{)}
\label{vanishing}
\label{tree}
Let $f\colon (X,C)\to (S,0)$ be a
Mori conic bundle. Then
\begin{enumerate}
\renewcommand\labelenumi{(\roman{enumi})}
\item  For an arbitrary ideal $\cc{I}$ such that
$\mt{Supp}(\cc{O}_X/\cc{I})\subset C$ we have,
$H^1(\cc{O}_X/\cc{I})=0$.
\item
The fiber $C$
is a tree of non-singular rational curves
(i.~e. $p_a(C')\le 0$ for any subcurve $C'\subset C$).
\item
If $C$ has $\rho$ irreducible components, then
$$
\mt{Pic}(X)\simeq H^2(C,\bb{Z})\simeq\bb{Z}^{\oplus\rho}.
$$
\qq
\end{enumerate}
\end{corollary1}

\subsection{Terminal singularities.}
\label{terminal}
Let $(X,P)$ be a terminal singularity of index $m\ge 1$
and let $\pi\colon (X\3,P\3 )\to (X,P)$ be the canonical cover.
Then $(X\3,P\3 )$ is a terminal singularity of index 1 and
$(X,P)\simeq (X\3,P\3)/\cyc{m}$.
It is known (see \cite{Pagoda}) that $(X\3,P\3)$ is a
hypersurface singularity, i.~e. there exist an
$\cyc{m}$-equivariant embedding
$(X\3,P\3)\subset (\CC^4,0)$.
Fix a character $\chi$ that generates $\mt{Hom}(\cyc{m},\CC^*)$.
For
$\cyc{m}$-semi-invariant  $z$ define weight
$\mt{wt}(z)$ as an integer defined  $\mod m$ such that
$$
\mt{wt}(z)\equiv a\mod m\qquad  {\rm iff}
\qquad
\zeta(z)=\chi(\zeta)^a\cdot z\quad
\text{for all}\quad \zeta\in\cyc{m}.
$$

\begin{theorem1}
\textup{(}see \cite{Mori-term}, \cite{YPG}\textup{)}
\label{cl-term}
In notations of \textup{\ref{terminal}}
 $(X\3,P\3)$ is analytically
$\cyc{m}$-isomorphic to
a hypersurface  $\phi=0$ in
$(\bb{C}^4_{x_1,x_2,x_3,x_4},0)$ such that $x_1,x_2,x_3,x_4$
are semi-invariants and we have one of the following two series:
\begin{itemize}
\item[{\rm (i)}]
\textup{(}the main series\textup{)}
$\mt{wt}(x_1,x_2,x_3,x_4;\phi )\equiv (a,-a,b,0;0)\mod m$, or
\item[{\rm (ii)}]
\textup{(}the exceptional case\textup{)}  ${m=4}$ and
$\mt{wt}(x_1,x_2,x_3,x_4;\phi )\equiv (a,-a,b,2;2)\mod 4$, \par\noindent
\end{itemize}
where $a$, $b$ are integers prime to $m$.
\end{theorem1}
Note that in case (i) the variety $X\3$ can be non-singular.
In this situation $(X,P)\simeq (\CC^3,0)/\cyc{m}(a,-a,b)$.
Such singularities are called {\it terminal cyclic
quotient singularities} and
denoted also by $\frac{1}{m}(a,-a,b)$.
\begin{theorem1}
\textup{(}see \cite{YPG}\textup{)}
\label{g.e.}
Let $(X,P)$ be a germ of terminal singularity.
Then a general member $F\in |-K_X|$ has only DuVal
singularity (at $P$).
\end{theorem1}
\subsubsection{}
\label{term-cl}
Terminal singularities can be classified in terms of
a general member
$F\in |-K_X|$ and its preimage $F\3\colon ={\pi}^{-1}(F)$
under
the canonical cover  $\pi\colon X\3\to X$ \cite{YPG}:
\par\bigskip
\begin{center}
{\large{
\begin{tabular}{|c|c|c|}
\hline
index&type&cover \quad $F\3\to F$\\
\hline
\multicolumn{3}{|c|}{{\rm the main series}}\\
\hline
$m$ & $cA/m $  & $A_{k-1}\stackrel{m:1}{\longrightarrow}A_{km-1} $\\
\hline
$2$ & $cAx/2$  & $A_{2k-1}\stackrel{2:1}{\longrightarrow}D_{k+2} $\\
\hline
$2$ & $cD/2 $  & $D_{k+1}\stackrel{2:1}{\longrightarrow}D_{2k}   $\\
\hline
$3$ & $cD/3 $  & $D_{4}\stackrel{3:1}{\longrightarrow}E_{6}      $\\
\hline
$2$ & $cE/2 $  & $E_{6}\stackrel{2:1}{\longrightarrow}E_{7}      $\\
\hline
\multicolumn{3}{|c|}{{\rm the exceptional case}}\\
\hline
$4$ & $cAx/4$  & $A_{2k-2}\stackrel{4:1}{\longrightarrow}D_{2k+1}$\\
\hline
\end{tabular}
}}
\end{center}
\par\bigskip
Cyclic quotient singularities are included in type $cA/m$.
\begin{definition1}
Let $X$ be a normal variety and $\mt{Cl}(X)$ be its
Weil divisor class group. The subgroup of
$\mt{Cl}(X)$ consisting of Weil divisor classes
which are $\QQ$-Cartier is called  the {\it semi-Cartier
divisor class group}. We denote it by $\mt{Cl}^{sc}(X)$.
\end{definition1}

\begin{theorem1}
\textup{(}see \cite{Pagoda},\cite{Kawamata}\textup{)}
\label{generator}
Let $(X,P)$ be a germ of 3-dimensional singularity.
Then $\mt{Cl}^{sc}(X,P)\simeq\cyc{m}$ and it is
generated by  the class of $K_{X}$.
Moreover the local fundamental group of
$(X,P)$ is cyclic:
$\pi_1(X-\{P\})\simeq\cyc{m}$.
\end{theorem1}
\begin{definition1}
Let $(S,0)$ be a two-dimensional log-terminal singularity.
It is well known (see \cite{Kawamata}) that
 $(S,0)$ is a quotient singularity,
i.~e. $(S,0)\simeq (\CC^2,0)/G$, where $G\subset GL(2,\CC)$ is a
finite subgroup without quasi-reflections.
The natural cover $(\CC^2,0)\to (S,0)$
we call the {\it topological cover} and the order of $G$ we call
the {\it topological index} of $(S,0)$.
\end{definition1}

The following is an easy consequence of \ref{terminal}.

\begin{lemma1}
\label{index-surf-threef}
 Let $(X,P)$ be a germ of a terminal threefold
singularity of index $m>1$ and
$(F,P)\subset (X,P)$ be a germ of
irreducible surface. Assume that $F$ is $\QQ$-Cartier and
$(F,P)$ is DuVal with topological index $n$.
\begin{itemize}
\item[{\rm (i)}]
Then $n$ is divisible by $m$.
\item[{\rm (ii)}]
 Moreover if $n=m$, then $(X,P)$ is a cyclic
quotient singularity and $(F,P)$ is of type $A_{m-1}$.
\end{itemize}
\end{lemma1}
\begin{proof}
Let $F\3:=\pi^{-1}(F)$. The divisor $F\3$ is  Cartier  on $X\3$, because
point $(X\3,P\3)$ is terminal of index 1
(see \cite[Lemma 5.1]{Kawamata}).
On the other hand $F\3$ is non-singular outside $P\3$.
Therefore $F\3$ is irreducible and normal. The cover $F\3\to F$ is \'etale
outside $P$. Hence the topological cover of $(F,P)$ can be facorized through
$F\3\to F$. This proves  (i). In conditions (ii) $F\3$ is non-singular
and so is $X\3$.
\end{proof}

\begin{proposition}
\textup{(}see  \cite{Mori-flip}\textup{)}
\label{cd}
Let $(X,P)$ be a germ of terminal singularity of
index $m$, $(C,P)\subset (X,P)$ be a germ of smooth curve and
 $\pi \colon (X\3 , P\3 )\to (X,P)$ be the canonical cover and
 $C\3:=(\pi^{-1}(C))_{\mt{red}}$. Then
\begin{itemize}
\item[{\rm (i)}]
 for  an arbitrary $\xi\in\mt{Cl}^{sc}(X,P)$,
there exists an effective (Weil) divisor $D$
such that $[D]=\xi$ and $\mt{Supp}(D)\cap C=\{ P\}$.
\item[{\rm (ii)}]
$\xi\to (D\cdot C)_P$ induces a homomorphism
$$
\mt{cl}(C,P)\colon  \mt{Cl}^{sc}(X,P)\to\frac{1}{m}\bb{Z}/\bb{Z}
\subset\bb{Q}/\bb{Z}.
$$
\end{itemize}
\qq
\end{proposition}
\begin{definition}
(see \cite{Mori-flip})
\label{shushu}
Let things be as in \ref{cd}. $X\supset C$ is called
{\it \textup{(}locally\textup{)} imprimitive} at $P$ of
{\it  splitting degree} $d$
if   $\mt{cl}(C,P)\colon \mt{Cl}^{sc}(X,P)\to\frac{1}{m}\bb{Z}/\bb{Z}$
is not an isomorphism and the order of
$\mt{Ker}(\mt{cl}(C,P))$ is equal $d$.
The order of $\mt{Coker}(\mt{cl}(C,P))$ is called the
{\it subindex} of $P$ and usually denoted by $\m$.
(Note that $\m d=m$).
 If   $\mt{cl}(C,P)\colon \mt{Cl}^{sc}(X,P)\to\frac{1}{m}\bb{Z}/\bb{Z}$
is an isomorphism, then $P$ is said to be
{\it \textup{(}locally\textup{)} primitive}.
In this case we put $d=1$, $\m=m$.
\end{definition}

In the situation above the preimage $C\3:=\pi^{-1}(C)$
under the canonical cover $\pi\colon X\3\to X$ has
exactly $d$ irreducible components.

\subsection{Construction, \cite{Mori-flip}.}
\label{Mori-covering}
Let  $f\colon (X,C\simeq\PP^1)\to (S,0)$ be a
 Mori conic bundle with irreducible central fiber and let
$P\in X$ be a point of index $m>1$.
 Assume that  $(X,C)$ is imprimitive of
splitting degree $d$ at $P$.
Take an effective Cartier divisor $H$
such that $H\cap C$ is a smooth point of $X$
and $(H\cdot C)=1$ and  take an effective Weil $\QQ$-Cartier divisor
$D$ such that $D\cap C=\{ P\}$ and
$D$ is a generator of $\mt{Cl}^{sc}(X,P)$ (see \textup{\ref{cd}}).  Then
the divisor
$$
D:=(m/d)D-((mD\cdot C)/d)H
$$
is a $d$-torsion element in $\mt{Cl}^{sc}(X,C)$. It defines a
finite Galois
$\cyc{d}$-morphism $g\6\colon X'\to X$ such that
$P':={g'}^{-1}(P)$
is one point, $g'$ is \'etale over $X-\{ P\}$
(hence $X'$ has only terminal singularities),
index of $(X',P')$ is equal to $m/d$, $C':=({g'}^{-1}(C))_{\mt{red}}$
is a union of $d$ $\bb{P}^1$'s meeting only at $P'$,
and each irreducible component of $C'$ is primitive at $P'$. \par

\section{Local invariants $w_P$ and $i_P$}
In this section we following Mori \cite{Mori-flip}
introduce numerical invariants $w_P$ and $i_P$.
There is nothing new in this section. The material is contained only
for convenience of the reader.
\subsection{}
Let $X$ be a normal three-dimensional
complex space with only terminal
singularities and let $C\subset X$ be a
reduced non-singular  curve. Denote by
$\cc{I}_C$ the ideal sheaf of $C$ and $\omega_X:=\OOO_X(K_X)$.
As in \cite{Mori-flip}, we consider the
following sheaves on $C$:
$$
\begin{array}{l}
\gr:=\hbox{{\rm torsion-free part of}}\ \omega_X/(\cc{I}_C\omega_X),\\
\\
\gro:=\hbox{{\rm torsion-free part of}}\ \cc{I}_C/\cc{I}^2_C.\\
\end{array}
$$
Since $C$ is non-singular, we have
$$
\omega_X/(\cc{I}_C\omega_X)=\gr\oplus\mt{Tors},
\qquad
\cc{I}_C/\cc{I}^2_C=\gro\oplus\mt{Tors}.
$$
\subsection{}
Let $m$ be the index of $X$.
The natural map
$$
(\omega_X\otimes\cc{O}_C)^{\otimes m}\to
\cc{O}_C(mK_X)
$$
induces an injection
$$
\beta\colon  (\gr)^{\otimes m}\to \cc{O}_C(mK_X).
$$
Denote
$$
w_P:=(\mt{length}_P\mt{Coker}\beta)/m.
$$
Note that if $f\colon (X,C\simeq\PP^1)\to (S,0)$ is a Mori conic bundle,
then $\deg \gr <0$, (because $\deg \cc{O}_C(mK_X)<0$).
\subsection{}
We also have the natural  map
$$
\begin{CD}
\gro\times\gro\times\omega_C
@>>>
\omega_X\otimes\cc{O}_C
\to
\gr,\\
x\times y\times zdt@>>> zdx\wedge dy\wedge dt\\
\end{CD}
$$
that induces a map
$$
\alpha\colon \wedge^2(\gro)\otimes
\omega_C\to\gr.
$$
Let
$$
i_P:=\mt{length}_P\mt{Coker}(\alpha ).
$$
Note that $i_P=0$ if $X$ is smooth at $P$.

\begin{lemma1}
\textup{(}\cite[2.15]{Mori-flip}\textup{)}
\label{non}
If $(X,P)$ is singular, then $i_P\ge 1$.\qq
\end{lemma1}

\begin{example}
\label{example}
Let $(X,P)$ be a cyclic quotient of type $1/\m(1,-1,\m+1)$
($\m$ is even) and let $\pi\colon (\CC^3_{x_1,x_2,x_3},0)\to (X,P)$ be
its canonical cover.
Consider the curve $C\3:=\{x_3=x_2^2-x_1^{2\m-2}=0\}\subset\CC^3$.
This curve has exactly two components
$C\3(1), C\3(2)=\{x_3=x_2\pm x_1^{\m-1}=0\}$ permuted under the action
of $\cyc{2\m}$.
Thus the image $C:=\pi(C\3)$ is a smooth irreducible curve,
$X$ is locally imprimitive along $C$ with splitting degree $2$.
The curve $C$ is naturally isomorphic to $C\3(1)/\cyc{\m}$
and a local uniformizing parameter on $C$ is $x_1^{\m}$.
It is easy to see that
\begin{itemize}
\item[(i)]
$\OOO_{C,P}\simeq\CC\{x_1^{\m}\}$.
\item[(ii)]
$\OOO_C(mK_X)\simeq
\OOO_C(d x_1\wedge d x_2\wedge d x_3)^m$,\qquad\qquad
$\gr \simeq \OOO_C(x_1^{\m-1}
x_1\wedge d x_2\wedge d x_3)^m$.
\item[(iii)]
$\gro\simeq\OOO_C(x_1^{\m-1}x_3)\oplus
\OOO_C(x_1^2(x_2^2-x_1^{2\m-2}))$.
\item[(iv)]
$w_P=(\m-1)/\m$,\quad  $i_p=2$.
\end{itemize}
\end{example}

From definitions we have

\begin{proposition}
\label{ocenki-1}
If $C\simeq\bb{P}^1$, then
$$
\deg\gro=2+\deg\gr -\sum_Pi_P,
$$
$$
(K_X\cdot C)=\deg\gr +\sum_Pw_P.\quad
$$
\qq
\end{proposition}
\begin{proposition}
\label{ocenki-2}
Let $f\colon (X,C\simeq\PP^1)\to (S,0)$ be a Mori conic bundle. Then
$$
\deg\gro\ge-2.
$$
\end{proposition}
\begin{proof}
Consider the exact sequence
$$
0\to\cc{I}_C/\cc{I}_C^{2}\to\cc{O}_X/\cc{I}_C^{2}
\to \cc{O}_C\to 0.
$$
By Corollary \ref{vanishing} $H^1(\cc{O}_X/\cc{I}_C^{2})=0$
and since $H^0(\cc{O}_X/\cc{I}_C^{2})
\to H^0(\cc{O}_C)$ is onto, we have
$H^1(\cc{I}_C/\cc{I}_C^{2})=0$.
Hence $H^1(\gro)=0$.
It gives us our assertion.
\end{proof}
\begin{corollary1}
\label{grw}
Let $f\colon (X,C\simeq\PP^1)\to (S,0)$ be a Mori conic bundle.
Then
$$
(-K_X\cdot C)+\sum w_P+\sum i_P\le 4.
$$
\end{corollary1}
\begin{remark3}
In the case of flipping contraction we always have
$\deg\gr=-1$ (see \cite{Mori-flip}) but in our situation
of Mori conic bundles
there are examples (e.~g. \ref{99} (i)) with $\deg\gr=-2$.
Nevertheless it is easy to see from \ref{ocenki-1} and \ref{ocenki-2}
that $-1\ge\deg\gr\ge -3$. J. Koll\'ar observe also that from vanishing
$R^1f_*\OOO_X(K_X+H)=0$ for any ample Cartier divisor $H$, we have
$\deg\gr =-1$ or $-2$.
\end{remark3}

\begin{corollary}
\label{<=3}
Let $f\colon (X,C\simeq\PP^1)\to (S,0)$ be a Mori conic bundle.
Then $(X,C)$
contains at most three singular points.
\end{corollary}
\begin{proof}
It follows from Corollary \ref{grw} and Lemma \ref{non}.
\end{proof}
As in \cite{Mori-flip} we need the following construction to
compute invariants $w_P$ and $i_P$ locally.
\subsection{}
Let $(X,P)$ be a terminal point of index $m$  and let
$C\subset (X,P)$ be a smooth curve.
Let $\pi\colon (X\3,P\3)\to (X,P)$ be the canonical cover and
$C\3:=(\pi^{-1}(C))_{\mt{red}}$.
We suppose that $(X,P)$ has splitting degree $d$ along $C$ and subindex
$\m$
(i.~e. $C\3$ has exactly $d$ components and $m=\m d$).
By \ref{terminal}, there exists an $\cyc{m}$-equivariant embedding
$X\3\subset\CC^4$. Let $\phi=0$ be an equation of $X\3$ in $\CC^4$.

 Let $C\3_1$ be an irreducible component of $C\3$.
For $z\in\OOO_{X\3}$ we define $\ord(z)$ as
the order of vanishing of $z$ on the normalization of $C\3_1$.
It is easy to see that $\ord(z)$ does not depend on our choice of
a component $C\3_1$,
if $z$ is a semi-invariant.
All the values  $\ord(z)$ form a semigroup which is denoted by
$\ord(C\3)$. Obviously the "coordinate" values $\ord(x_1),\dots,\ord(x_4)$
generate $\ord(C\3)$.

\begin{proposition-definition1}
\label{notations}
\label{nor}
By \cite{Mori-flip}, one can
choose a coordinate system $x_1, x_2, x_3, x_4$ in
$\CC^4$ and a character $\chi\colon \cyc{m}\to\CC^*$
such that the following conditions hold.
\begin{itemize}
\item[{\rm (i)}]
Coordinates $x_1, x_2, x_3, x_4$ are semi-invariants with
weights satisfying \textup{\ref{cl-term}}.
\item[{\rm (ii)}]
$\ord(x_i)\equiv\wt(x_i)\mod \m$ for all $i=1,2,3,4$.
\item[{\rm (iii)}]
Each component of $C\3$ is parameterized as
$$
x_i= \chi(g)^{\wt(x_i)}t^{\ord(x_i)},\quad i=1,2,3,4,\qquad
g\in\cyc{m}/\cyc{\m}.
$$
\item[{\rm (iv)}]\label{normal}
{\it Normalizedness property.}
For any $i=1,2,3,4$ there is no semi-invariants $y$
such that $\wt(y)\equiv\wt(x_i)\mod m$ and $\ord(y)<\ord(x_i)$.
In particular $\ord(x_i)<\infty$ for all $i$.
\item[{\rm (v)}]
There exists an invariant function $z$ on $X\3$ such that
$\ord(z)=\m$.
In particular in  the case of the main series one has
$\ord(x_4)=\m$, $\wt(x_4)\equiv 0\mod m$.
 \end{itemize}
Such a coordinate system is said to be {\it normalized}.
Note that we still may  permute $x_1$, $x_2$ and may replace
a character $\chi$ with $\chi'$ if $\chi'\equiv \chi\mod \m$.\qq
\end{proposition-definition1}

\begin{proposition}
\textup{(}see \cite[2.10]{Mori-flip}\textup{)}
\label{computation-w}
Under the notations and conditions \textup{\ref{notations}}
 one has
$$
w_P=\min \{ \ord(\psi)\quad |\quad  \psi\in\OOO_{X\3,P\3},\quad
\wt(\psi)\equiv-\wt(x_3)\mod m\}.
$$
\qq
\end{proposition}

\begin{proposition}
 \textup{(}cf. \cite[0.4.14.2]{Mori-flip}\textup{)}
\label{computation-k}
Let $F\in |-K_{(X,P)}|$ be a general member. Then in a normalized
coordinate system we have
$(F\cdot C)_P=\ord(x_3)/\m$.
\qq
\end{proposition}

\begin{proposition}
\textup{(}see \cite[2.12]{Mori-flip}\textup{)}
\label{computation-i}
Under the notations of \textup{\ref{notations}},
one has
$$
i_P\m=\m-\ord(x_4)-\m w_P+\min_{\phi_1,\phi_2\in
J\3_0}[\phi,\phi_1,\phi_2],
$$
where $J\3_0$ is the invariant part of the ideal sheaf
$J\3$ of $C\3$ in $\CC^4$, $\phi$ is an equation of
$X\3$ in $\CC^4$ and
$$
[\psi_1,\psi_2,\psi_3]:=
\ord\partial (\psi_1,\psi_2,\psi_3)/\partial (x_1,x_2,x_3).
$$
is the Jacobian determinant.\qq
\end{proposition}
By semi-additivity of $[\ ,\ ,\ ]$ and
because $x_4\phi\in J\3_0$, we have

\begin{corollary1}
\textup{(}see \cite[(2.15.1)]{Mori-flip}\textup{)}
\label{predv}
In notations above
$$
i_P\m\ge -\m w_P+\min_{\psi_1,\psi_2,\psi_3
\in J\3_0}[\phi_1,\phi_2,\phi_3]. \quad
$$
\qq
\end{corollary1}
\subsubsection{}
\label{choice}
It is clear that the ideal $J\3_0$ is generated by
invariants of the form
$\psi-x_4^n$ where $\psi$ is an invariant  monomial with
$\ord(\psi)=\ord(x_4^n)=n\m$. Thus we can assume in \ref{predv}
that $\phi_i=\psi_i-x_4^{n_i}$, for $i=1,2,3$, where
each $\psi_i$ is an invariant monomial in $x_1,x_2,x_3$.
For any element
$\psi=x_1^{b_1}x_2^{b_2}x_3^{b_3}x_4^{b_4}-x_4^{n}$ define
the vector $\ex(\psi)\in\ZZ^3$ by $\ex(\psi)=(b_1, b_2, b_3)$.

\begin{lemma1}
\textup{(}see \cite[2.14]{Mori-flip}\textup{)}
Let $\psi_1,\psi_2,\psi_3$ be as in \textup{\ref{choice}}.
Then
$$
[\psi_1,\psi_2,\psi_3]=
\left\{
\begin{array}{ll}
\sum_{i=1}^3\ord(\psi_i)-\sum_{i=1}^3\ord(x_i)
 & \hbox{{\rm if}}\quad \ex(\psi_i), i=1,2,3 \\
 &\hbox{{\rm  are independent,}}\\
 &\\
\infty     &\hbox{{\rm otherwise}}.\quad
\end{array}
\right.
$$
\qq
\end{lemma1}

An invariant monomial $\psi$ in $x_1, x_2, x_3$ is said
to be {\it simple}
if it cannot be presented as a product of two non-constant
invariant monomials.

\begin{corollary1}
\textup{(}see \cite[proof of 2.15]{Mori-flip}\textup{)}
\label{formula}
Under the notations and conditions of \textup{\ref{notations}}
for some simple different invariant monomials $\psi_1$,
$\psi_2$, $\psi_3$
 in $x_1, x_2, x_3$ we have
$$
\begin{array}{l}
\m i_P\ge(\ord(\psi_1)-\ord(x_1))+(\ord(\psi_2)-\ord(x_2))+\\
\qquad (\ord(\psi_3)-\ord(x_3)-\m w_P).
\end{array}
\leqno(*)
$$
 Moreover up to permutations $\psi_1$, $\psi_2$, $\psi_3$
we may assume that $\psi_i=x_i\nu_i$, where $\nu_i$, $i=1,2,3$ also
are  monomials. In particular,
all three terms in the formula are non-negative and
$$
\m i_P\ge\ord(\nu_1)+\ord(\nu_2)+\ord(\nu_3)-\m w_P.
\leqno(**)
$$
\qq
\end{corollary1}

\section{Easy lemmas}
\subsection{Construction.}
\label{covering}
Let  $f\colon (X,C)\to (S,0)$ be a
 Mori conic bundle.
Assume that $(S,0)$ is singular.
 Then $(S,0)$ is a log-terminal point
(see e.~g. \cite{KoMM} or \cite{Ishii}) and
the topological cover $h\colon (S\6,0)\simeq
(\bb{C}^2,0)\to (S,0)$ is non-trivial.
Let $X\6$ be the normalization of $X\times_{S}S\6$ and
 $G=\mt{Gal}(S\6/S)$. Then we have the diagram
$$
\begin{CD}
X\6@>{g}>>X\\
@V{f\6}VV @V{f}VV\\
S\6@>{h}>>S\\
\end{CD}
$$
The group $G$ acts on $X\6$ and clearly
$X=X\6/G$. Since the action of $G$ on $S\6-\{ 0\}$ is free,
so is the action of $G$ on $X\6-C\6$, where
$C\6:=({f\6}^{-1}(0))_{\mt{red}}$.
Therefore  $X\6$ has only terminal singularities
and the induced action of $G$ on $X\6$ is free
outside of a finite set
of points (see e.~g. \cite[6.7]{CKM}).
Since $K_{X\6}=g^*(K_X)$,
$f\6\colon (X\6,C\6)\to (S\6,0)$ is a Mori conic bundle
 with non-singular base surface.
\par
As an easy consequence of the construction above and \ref{generator}
we have
\begin{proposition}
\textup{(}see \cite{Pro}, \cite{Pro1}, \cite{Kollar}\textup{)}
\label{ooo}
Let  $f\colon (X,C)\to (S,0)$   be a  Mori conic bundle.
Then $(S,0)$ is a cyclic quotient singularity.
\end{proposition}

\begin{corollary}
Let  $f\colon (X,C)\to (S,0)$   be a  Mori conic bundle
and let   $T$ be the torsion part of the Weil divisor class group
$\mt{Cl}(X)$. Then
\begin{itemize}
\item[{\rm (i)}]
$T$  is a cyclic group,
\item[{\rm (ii)}]
the topological index of $(S,0)$ is equal to the order of $T$.
\end{itemize}
In particular, $(S,0)$ is non-singular iff $T=0$.
\end{corollary}
\begin{proof}
Consider the diagram from  \ref{covering}. We claim that $\mt{Cl}(X\6)$
is torsion-free. Indeed in the opposite case there is
a finite \'etale in codimension 1 cover $X'\to X\6$ and  then it must
be \'etale outside finite number
singular points on $C\6$. One can take the Stein factorization
$X'\to S'\to S\6$. Then $S'\to S\6$ will be \'etale outside $0$.
This contradicts smoothness of $S\6$.  Therefore $\mt{Cl}(X\6)$
is torsion-free and $X\6\to X$ is the maximal abelian cover of
$X$. It gives us $T\simeq\cyc{n}$, where $n$ is the topological index of
$(S,0)$.
\end{proof}

\begin{proposition}
Let  $f\colon (X,C\simeq\PP^1)\to (S,0)$ be a
 Mori conic bundle with irreducible central fiber.
Suppose that the point $(S,0)$ is singular of topological index $n>1$.
Then for $f\6 \colon (X\6,C\6)\to (S\6,0)$ of \textup{\ref{covering}}
we have the following two possibilities.
\begin{itemize}
\item[{\rm (i)}] \textup{(}primitive case\textup{)}
$C\6$ is
irreducible, $X$ has at least two non-Gorenstein points.
\item[{\rm (ii)}] \textup{(}imprimitive case\textup{)}
$C\6$ has $d>1$
irreducible components, they all pass through one point $P\6$ and
they do not intersect elsewhere.  If in this case $m\6$ is the index
of $P\6$ and $n$ is the topological index of $(S,0)$, then  the index
of the point $g(P)\in X$ is equal to $nm\6$.  \end{itemize}
\end{proposition} \begin{proof} If $C\6$ is irreducible, then  the
action of $G\simeq\cyc{n}$ on $C\6\simeq\PP^1$ has exactly two fixed
points $P_1\6$ and  $P_2\6$. Their images $g(P_1\6)$ and  $g(P_2\6)$
are points of indices $>1$. Assume that $C\6$ is reducible.  Then
$G\simeq\cyc{n}$  acts on components $\{C\6_i\}$ of $C\6$
transitively.  Since $\cup C\6_i$ is a tree, $\cap C\6_i$ is a point
and this point must be fixed under $G$-action. The rest is obvious.
\end{proof}

\begin{corollary}
\label{odnatochka}
Let  $f\colon (X,C\simeq\PP^1)\to (S,0)$ be a
 Mori conic bundle with irreducible central fiber.
Assume that $X$ contains only one non-Gorenstein point $P$.
Then the following conditions are equivalent
\begin{itemize}
\item[{\rm (i)}]
$(S,0)$ is singular of topological index $n$,
\item[{\rm (ii)}]
$C$ is locally imprimitive
\textup{(}see
\textup{\ref{shushu}}\textup{)}
at $P$ of splitting degree $n$.
\end{itemize}
\end{corollary}

\begin{lemma}
\label{s-even}
\label{even}
Let $f\colon (X,C\simeq\PP^1)\to (S,0)$ be a minimal Mori conic bundle.
Assume that $f$ is imprimitive of splitting degree $d$. Then $d$ is
even.
\end{lemma}
\begin{proof}
Consider the base change \ref{covering}.
  Then $f\6\colon (X\6,C\6)\to (S\6,0)$ is a Mori conic
  bundle with  $\rho (X\6)=d$.
 Since $-K_{X\6}$ is $f\6$-ample, the Mori cone
$\overline{NE}((X\6,C\6)/(S\6,0))\subset \bb{R}^{\rho}$
is generated by classes of $C_1\6,\dots, C_d\6$. Thus each $C_i\6$
generates an extremal ray $R_i$ and
$\overline{NE}((X\6,C\6)/(S\6,0))\subset \bb{R}^{\rho}$
is simplicial. The  contraction of each extremal face of
$\overline{NE}((X,C)/(S,0))$
gives us  an extremal neighborhood in the sence of
\cite{Mori-flip} (not necessary flipping).
As in the introduction we
apply the Minimal Model Program to $(X\6,C\6)$ over $(S\6,0)$:
$$
\begin{array}{ccccc}
(X\6,C\6)&         &\stackrel{p}{\bir}&        &(Y,L)\\
         &\searrow&              &\swarrow&     \\
         &         &(S\6,0)       &        &     \\
\end{array}
$$
Here $(Y,L)\to (S\6,0)$ is a Mori conic bundle with $\rho(Y,S\6)=1$.
Then the map $p$ contracts at least $d-1$ divisors.
Let $E$ be a proper transforms one of them on $X\6$.
Since $f\6$ is flat $\Gamma :=f\6(E)$ is a curve on $S\6$.
Moreover for a general point $s\in\Gamma$ the preimage
${f\6}^{-1}(s)$ is a reducible conic, so ${f\6}^{-1}(s)=\ell_1+\ell_2$.
Therefore ${f\6}^{-1}(\Gamma)$ has exactly two irreducible components
$E\supset\ell_1$
and $E'\supset\ell_2$.
Consider the orbit $\{E=E_1, E_2, \dots, E_k\}$ of $E$ under
the action of $\cyc{d}$.
Obviously, every $f\6(E_i)$ is a curve on $S\6$.
We have $\sum E_i\sim aK_{X\6}$ for some $a\in\QQ$
because $\cyc{d}$-invariant part of
$\mt{Pic}(X\6)\otimes\QQ$ is $\mt{Pic}(X)\otimes\QQ\simeq\QQ$.
But for a general fiber $\ell:={f\6}^{-1}(s)$, over  $s\in S\6$ we have
$(\ell\cdot\sum E_i)=0$, $(\ell\cdot K_{X\6})<0$, so $a=0$,
and $\sum E_i\sim 0$. Further $(\ell_1\cdot E)<0$, $(\ell_1\cdot E')>0$.
It gives us  $(\ell_1\cdot E_i)>0$ for some
$E_i\in\{E=E_1, E_2, \dots, E_k\}$. Then $E'=E_i$, i.~e.
there exists $\sigma\in\cyc{d}$ such that $\sigma(E)=E'$.
From the symmetry we get that the orbit $\{E=E_1, E_2, \dots, E_k\}$
may be divided into couples of divisors $E_j$, $E_j'$ such that
${f\6}(E_j)={f\6}(E_j')$ is a curve. Thus both $k$ and $d$ are even.
\end{proof}

\begin{lemma}
\label{warwick}
Let $f\colon (X,C\simeq\PP^1)\to (S,0)$ be a Mori conic bundle
with irreducible  central fiber and with a unique point $P$ of index
$m>1$. Assume that   $P$ is an imprimitive of subindex $\bar{m}$
and splitting degree $d$ \textup{(}so $m=\bar{m}d$\textup{)}.
Then
\begin{itemize}
\item[{\rm (i)}]
$2\bar{m}\equiv 0\mt{mod} d$,
\item[{\rm (ii)}]
$(-K_X\cdot C)=1/\m$.
\end{itemize}
\end{lemma}
\begin{corollary1}
\label{warwick1}
In conditions \textup{\ref{warwick}}, $\m\ge d/2$.
In particular, if $d\ge 3$, then $\m>1$.
\end{corollary1}
\begin{proof}
\label{xaxa}
Again consider the base change \ref{covering}.
 Terminal singularities are rational
and therefore Cohen-Macaulay \cite{Kempf}.
By  \cite[Theorem 23.1]{Matsumura}, $f\6$ is flat.
Let $X\6_s$ be a general fiber of $f\6$.
Then $X\6_s\equiv nC\6$ for some $n\in \NN$.
We have
$$
2=(-K_{X\6}\cdot X\6_s)=n(-K_{X\6}\cdot C\6)=n(-K_{X\6}\cdot \sum C\6(i))=
nd(-K_{X\6}\cdot C\6(1)).
$$
It is clear that
 $(-K_{X\6}\cdot  C\6(1))\in\frac{1}{\m}\ZZ$ because $-\m K_{X\6}$
is Cartier.
If $(-K_{X\6}\cdot  C\6(1))=\delta /\bar{m}$ for some $\delta \in\NN$,
then $2=nd\delta /\bar{m}$ or equivalently $2\bar{m}=nd\delta $.
This proves (i).
Since every $(X\6,C(i))$ is a primitive  extremal neighborhood
and $-K_{X\6}$ is a generator of
$\mt{Cl}^{sc}(X\6,P\6)\simeq\cyc{\bar{m}}$, $(\delta ,\bar{m})=1$.
Thus $\delta=1$ by Lemma \ref{even}.
 On the other hand,
$(-K_{X}\cdot  C)=(-K_{X\6}\cdot  C\6(1))=\delta /\bar{m}$,
because $K_{X\6}=g^*(-K_{X})$.
The rest is obvious.
\end{proof}

Using Lemma \ref{computation-k} we obtain.
\begin{corollary1}
\label{a3p}
Let $f\colon (X,C\simeq\PP^1)\to (S,0)$ be a Mori conic bundle.
Assume that  $(S,0)$ is singular and
$X$ contains the only one non-Gorenstein
point $P$. Then in notations of \textup{\ref{nor}} one has
$\ord(x_3)\equiv 1\mod \m$.\qq
\end{corollary1}

The following construction explains why Conjecture $(*_2)$ can
be useful for studying Mori conic bundles.

\subsection{Construction (double cover trick \cite{Kawamata}).}
\label{construction-II}
Let $f\colon (X,C)\to (S,0)$ be a Mori conic bundle and let
$D\in |-2K_X|$ be a general member.
Assume that $K_X+1/2D$ is log-terminal (i.~e. $(*_2)$ holds).
Then there exists a double cover $h\colon Y\to X$ with
ramification divisor
$D$ such that $Y$ has only canonical Gorenstein singularities
(see \cite[8.5]{Kawamata}).
Thus we have two morphisms
$$
Y\stackrel{h}{\longrightarrow}X
\stackrel{f}{\longrightarrow} S.
$$
By the Hurwitz formula, $K_Y=h^*(K_X+1/2D)=0$.
Therefore the composition map $g\colon Y\to S$ gives us
an elliptic fiber space structure on $Y$.

\section{The main result}
The main result of this section is the following
\begin{theorem}
\label{main}
Let $f\colon (X,C\simeq\PP^1)\to (S,0)$ be a
 Mori conic bundle with only one non-Gorenstein point $P$.
Then one of the following holds
\subsubsection{}
\label{main0}
$(S,0)$ is non-singular.

\subsubsection{}
\label{main1}
 $(S,0)$ is a DuVal point of type $A_1$ and splitting degree
of $P$ is $2$ and $X$ contains at most one more singular
\textup{(}Gorenstein\textup{)} point.
There are the following subcases:
  \begin{itemize}
  \item[{\rm (i)}]
 \textup{(}cf. \cite{Pro2}, \cite{Pro3}\textup{)}
\label{war2}\label{99}
$f\colon (X,C)\to (S,0)$ is a quotient of
a conic bundle $f\6\colon (X\6,C\6)\to (\CC^2,0)$ by
$\cyc{2}$, $P$ is
the only singular point.
In some coordinate system $(x,y,z;u,v)$ in $\PP^2\times\CC^2$
the variety $X\6$ can be defined by the equation
$$
x^2+y^2+\psi(u,v)z^2=0, \qquad \psi(u,v)\in\CC\{u,v\},
$$
where $\psi(u,v)$ has no multiple factors
and  consists only of monomials of even degree.
Here
$X\6\to\CC^2$ is the natural projection and the
action of $\cyc{2}$ has the form
$$
u\to -u,\qquad v\to -v,\qquad x\to -x,\qquad y\to y,\qquad z\to z.
$$
Moreover
\par
\mbox{{\rm (a)}}
 if $\mt{mult}_{(0,0)}\psi(u,v)=2$, then $(X,P)$ is of type $cA/2$,
\par
\mbox{{\rm (b)}}
if $\mt{mult}_{(0,0)}\psi (u,v)\ge 4$, then $(X,P)$ is of type $cAx/2$.

\item[{\rm (ii)}] $(X,P)$ is of type $cAx/4$, more precisely
$$
(X,P)\simeq
(\{y_1^2-y_2^2+y_3\varphi_3+y_4\varphi_4=0\}/\cyc{4}(1,3,3,2),
$$
$$
C= \{y_1^2-y_2^2=y_3=y_4=0\})/\cyc{4},
$$
where $\varphi_3$, $\varphi_4$ are semi-invariants with weights
$$
\wt(\varphi_3)\equiv 3\mod 4,\qquad
\wt(\varphi_4)\equiv 0\mod 4
 $$
\textup{(}this is a point of type $II^{\vee}$ in
Mori's classification \cite{Mori-flip}
\textup{)}.
\item[{\rm (iii)}]
$(X,P)$ is the cyclic quotient singularity
$\CC^3/\cyc{4}(1,3,1)$ and
$$
C=\{y_1^2-y_2^2=y_3=0\}/\cyc{4}.
$$
\item[{\rm (iv)}] $\m$ is even  and $\ge 4$,
$(X,P)$ is the cyclic quotient singularity
$\CC^3/\cyc{2\m}(1,-1,\m+1)$
$$
C=  \{ y_2^2-y_1^{2\m-2}=y_3=0\}\ni 0)/\cyc{2\m}
$$
\textup{(}this is a point of type $IC^{\vee}$ in
\cite{Mori-flip}, see also example
\textup{\ref{example}}\textup{)}.
\item[{\rm (v)}] $\m$ is even, $(X,P)$ is of type $cA/m$,
where $m=2\m$. More precisely
$$
(X,P)\simeq \{y_1y_2+y_2\varphi_1+y_3\varphi_2
+(y_4^2-y_1^{2\m})\varphi_3
=0\}/\cyc{2\m}(1,-1,\m+1,0)
$$
 and
$$
C=\{y_2=y_3=y_4^2-y_1^{2\m}=0\}/\cyc{2\m},
$$
where $\varphi_i$, $i=1,2,3$ are semi-invariants with weights
$$
\wt(\varphi_1,\varphi_2,\varphi_3)\equiv (1,\m -1,0) \mod 2\m.
$$
 \end{itemize}
\subsubsection{}
\label{main2}
 $(S,0)$ is a DuVal point of type $A_3$,
$(X,P)$ is a cyclic quotient singularity of type
$\CC^3/\cyc{8}(a,-a,1)$,
splitting degree of $P$ is $4$
and $X$ has no other singular points. There are two subcases
\begin{itemize}
  \item[{\rm (i)}]
 $a=1,\qquad\qquad C=\{y_1^{4}-y_2^{4}=y_3=0\}/\cyc{8}$,
 \item[{\rm (ii)}]
$a=3,\qquad\qquad C=\{y_2^2-y_3^2=y_2y_3-y_1^2=0\}/\cyc{8}$.
  \end{itemize}

\end{theorem}

\begin{remark}
We will show that
in  \ref{main1}, \ref{main2} a general member of
$|-K_X|$ does not contain $C$ except \ref{main1} (iv), (v).
In these two cases Conjecture $(*_2)$ holds.
Using the same arguments one can prove that in \ref{main0}
at least one of Conjectures $(*_1)$, $(*_2)$ or $(*_3)$ holds.
More precise results for the case of non-singular base surface
will be published elsewhere.
\end{remark}
Before starting the proof we consider examples of Mori conic bundles
such as in \ref{main1} (ii), (iii) and \ref{main2}.
Unfortunately I don't know  examples of Mori conic bundles
as in \ref{main1} (iv) or (v).

\begin{example}
Let $Y\subset\PP^3_{x,y,z,t}\times\CC^2_{u,v}$ be a non-singular
subvariety given by the equations
$$
\left\{
\begin{array}{l}
xy=ut^2\\
(x+y+z)z=vt^2
\end{array}
\right.
$$
It is easy to check that $Y$ is non-singular.
The projection $g\colon Y\to\CC^2$ gives us an elliptic fibration
whose  fibers are intersections of two quadrics in $\PP^3$. The
special fiber $g^{-1}(0)$ is the union of four lines meeting at one point
$Q:=\{(x,y,z,t;u,v)=(0,0,0,1;0,0)\}$.
Consider also the following action of $\cyc{8}$ on $Y$:
\begin{multline*}
x\to\varepsilon^{-3}z,\qquad y\to\varepsilon(x+y+z),
\qquad  z\to-\varepsilon y,\qquad
t\to t,\\
u\to \varepsilon^{-2}v,\qquad  v\to-\varepsilon^{2}u,
\end{multline*}
where $\varepsilon:=\exp(2\pi i/8)$.
Then a generator of $\cyc{8}$ has only one fixed point $Q$ and the quotient
$(Y,Q)$ is terminal of type $\frac{1}{8}(1,-1,3)$.
The locus of fixed point for $\cyc{2}\subset\cyc{8}$
consists of the point $Q$ and the divisor $B:=\{t=0\}\cap Y$.
Let $h\colon Y\to Y':=Y/\cyc{2}$ be the quotient morphism.
By the Hurwitz formula $0=K_Y=h^*K_{Y'}+B$, hence $h^*(-K_{Y'})=B$.
Therefore $-K_{Y'}$ is ample over $\CC^2$ and $Y'\to\CC^2$ is a Mori conic
bundle with singular point $h(Q)$ of type $\frac{1}{2}(1,1,1)$ and
reducible central fiber. Note that
$$
g\colon Y\stackrel{h}{\longrightarrow}X
\stackrel{f}{\longrightarrow} S.
$$
is the double cover trick \ref{construction-II}.
Further the induced action $\cyc{4}=\cyc{8}/\cyc{2}$ on $Y'$
is free outside point $h(Q)$ and permutes components of the central fiber.
Therefore the anticanonical divisor of $X:=Y'/\cyc{4}$ is ample over
$\CC^2/\cyc{4}$ and the central fiber of the projection
$X\to\CC^2/\cyc{4}$ is irreducible.
Finally we obtain the  Mori conic bundle  $X:=Y/\cyc{8}\to\CC^2/\cyc{4}$
such as in \ref{main2}.
The diagram
$$
\begin{CD}
Y'@>>>X\\
@V{f'}VV @V{f}VV\\
\CC^2@>>>\CC^2/\cyc{4}
\end{CD}
$$
is exactly as in  the construction  \ref{covering}.
\end{example}

\begin{example}
Let  $Y\subset\PP^3_{x,y,z,t}\times\CC^2_{u,v}$
is given by the equations
$$
\left\{
\begin{array}{l}
xy=ut^2\\
z^2=u(x^2+y^2)+vt^2
\end{array}
\right.
$$
 Consider the action of  $\cyc{4}$ on $Y$:
$$
x\to y,\qquad y\to -x,
\qquad  z\to iz,\qquad
t\to t,\qquad  u\to-u,\qquad  v\to-v.
$$
As above
 $X:=Y/\cyc{4}\to\CC^2/\cyc{2}$
is a Mori conic bundle such as in \ref{main1} (iii)
with unique singular point
of type
 $\frac{1}{4}(1,1,3)$
\end{example}

\begin{example}
Let  $Y\subset\PP^3_{x,y,z,t}\times\CC^2_{u,v}$
is given by the equations
$$
\left\{
\begin{array}{l}
xy=(u^{2k+1}+v)t^2,\qquad k\in\NN\\
z^2=u(x^2-y^2)+vt^2
\end{array}
\right.
$$
Define the action of  $\cyc{4}$ on $Y$ by
$$
x\to iy,\qquad y\to ix,
\qquad  z\to iz,\qquad
t\to t,\qquad  u\to-u,\qquad  v\to-v.
$$
Then
 $X:=Y/\cyc{4}\to\CC^2/\cyc{2}$
is a Mori conic bundle such as in \ref{main1} (ii)
with unique singular point
of type
 $cAx/4$.
\end{example}
\subsection{}
To begin the proof of Theorem \ref{main} we
suppose that $(S,0)$ is singular
(otherwise we have case \ref{main0}). Then by \ref{odnatochka}
the point $P$ is imprimitive.
Denote by $m$, $d$ and  $\m$ its index, splitting degree and  subindex,
respectively.
First we prove the existence of good divisors in $|-K_X|$
(with two exceptions \ref{main1} (iv), (v)).
As in \cite[\S 7]{Mori-flip} we will use local methods to find
good divisors in $|-K_X|$ or  $|-2K_X|$.
\begin{lemma}
\textup{(}cf. \cite[ proof of 7.3]{Mori-flip}\textup{)}.
\label{assumptions}
Let $f\colon (X,C\simeq\PP^1)\to (S,0)$ be a
 Mori conic bundle with only one non-Gorenstein point $P$.
Let $\pi\colon  (X\3,P\3)\to (X,P)$ be the canonical
cover and $X\3\subset\CC^4$
 be an embedding  as in \textup{\ref{cl-term}}.
Then a general member $F\in |-K_X|$ does not contain $C$
and has  only DuVal singularity at $P$ iff
we have $\ord(x_3)<\m$ in a normalized coordinate system
$(x_1,\dots,x_4)$ in $\CC^4$. \textup{(}By \textup{\ref{a3p}},
$\ord(x_3)<\m$
implies $\ord(x_3)=1$, if $(S,0)$ is singular.\textup{)}
\end{lemma}
\begin{proof}
Assume that $\ord(x_3)<\m$.
 From Lemma \ref{warwick} we have $(-K_X\cdot C)=1/\m$.
Take a general member $F\in|-K_{(X,P)}|$. Proposition
\ref{computation-k} yields $(F\cdot C)_P=\ord(x_3)/\m<1$.
Then $F+K_X$ can be considered as
a Cartier divisor on $X$ and in our  case
 $-1<((F+K_X)\cdot C)<1$.  Therefore $(F\cdot C)=(-K_X\cdot C)$.
Since $\mt{Pic}(X)\simeq\ZZ$, $F\in|-K_X|$  and by Theorem
 \ref{g.e.}, $F$ has only DuVal singularity at $P$.
The inverse implication is obvious.
\end{proof}
Thus it is sufficient to show only  $\ord(x_3)<\m$.
First we consider the exceptional series of terminal points
\ref{cl-term}, (ii).
\begin{lemma}
\textup{(}cf. \cite[4.2]{Mori-flip}\textup{)}.
\label{except}
If $(X,P)\supset C$ is an imprimitive point of type $cAx/4$, then
 $\m=2$ and
$$
\begin{array}{cccccc}
              &x_1&x_2&x_3&x_4&\\
\mt{wt}  &1  &3  &3  &2 &\mod 4 \\
\mt{ord}  &1  &1  &1  &2 & \\
\end{array}
$$
\end{lemma}
\begin{proof} (cf. \cite[3.8]{Mori-flip})
From Corollary \ref{warwick1}, one has $\m=2$.
By \ref{notations}, $\ord(x_4)\equiv 0\mod \m$ and
there exists an invariant monomial $\psi$ such that $\ord(\psi)=\m=2$.
Since $\wt(x_i)\not\equiv 0\mod 4$ for $i=1,2,3,4$, we have $\psi=x_ix_j$.
Up to permutation $\{1, 2\}$, we may assume that
$\wt(x_2)\equiv\wt(x_3)$. Therefore $\ord(x_2)=\ord(x_3)$.
It gives us $\psi=x_1x_2$ or $x_1x_3$, whence
$\ord(x_1)=\ord(x_2)=\ord(x_3)=1$.
If necessary we can replace the character $\chi$ with
$-\chi$ to obtain $\wt(x_1)\equiv 1$.
Finally
$\wt(x_4)\equiv\wt(x_1^2)\equiv 2\mod m$ gives as $\ord(x_4)=2$.
\end{proof}
Since $\ord(x_3)=1$, we have a good divisor in $|-K_X|$ in this case.
By \ref{nor}, components $C\3(1), C\3(2)$ of curve $C\3$ are
images of
$$
t\longrightarrow (t,t,t,t^2),\qquad t\longrightarrow (it,-it,-it,-t^2).
$$
Therefore $C\3$ can be given by the equations
$x_4-x_1^2=x_2-x_3=x_2^2-x_1^2=0$. By changing coordinates as
$y_1=x_1$, $y_2=x_2$,
$y_3=x_3-x^2$,
$y_4=x_4-x_1^2$ we obtain case \ref{main1} (ii) of our theorem.
\par
From now we assume that $(X,P)$ is from the main series.
In particular $\wt(x_4)\equiv 0\mod m$, $\ord(x_4)=\m$.

\begin{lemma}
\label{matrix}
Let $f(X,C\simeq \PP^1)\to (S,0)$ be a Mori conic bundle.
Assume that $P\in X$ is the only point of index $m>1$.
Assume further that $(X,C)$ is imprimitive at $P$ of splitting degree $d$.
Then one of the following holds
\subsubsection{}
\label{matrix1}
A general member of $|-K_X|$ does not contain $C$ and
has DuVal singularity at $P$.
\textup{(}By \textup{\ref{assumptions}}, this is equivalent to
$\ord(x_3)=1$\textup{)}.
\subsubsection{}
\label{matrix2}
$d=2$, $\m$ is even
 a general member of $|-2K_X|$ does not contain $C$,
and up to permutation $x_1,\dots,x_4$ in some
normalized coordinate system we have
\begin{itemize}
\item[{\rm (i)}]
\quad
$
\begin{array}{llllll}
      & x_1   & x_2   & x_3   & x_4& \\
\ord  & 1     & \m-1  & \m+1  & \m & \\
\wt   & 1     & -1 & \m+1  & 0   &\mod m\\
\end{array}
$
\qquad\qquad  $\m\ge 4$\quad or
\item[{\rm (ii)}]
\quad
$
\begin{array}{llllll}
      & x_1   & x_2   & x_3   & x_4 &\\
\ord  & 1     & 2\m-1 & \m+1  & \m  &\\
\wt   & 1     & -1 & \m+1  & 0  &\mod m \\
\end{array}
$
\end{itemize}
\end{lemma}

\subsection{Proof.}
\label{assumptions1}
Denote $a_i:=\ord(x_i)$. Since $(\wt(x_i), m)=1$, for $i=1,2,3$,
one has $(a_i,\m)=1$, $i=1,2,3$.
For \ref{matrix2} it is sufficient to prove that
$a_3<\m$. Let $a_3>\m$.
First suppose that $\ord(x_1^m)=2\m$.
We claim that in this situation  \ref{matrix2} holds.
Indeed, then
$a_1m=a_1\m d=2\m$, hence $d=2$, $a_1=1$, $m=2\m$.
Thus $\wt(x_1)\equiv\wt(x_3)\equiv 1\mod\m$.
By the assumptions in \ref{normal} and because $a_3>\m$,
$\wt(x_1)\not\equiv
\wt(x_3)\mod m$. Since $(\wt(x_1), m)=(\wt(x_3), m)=1$,
$\m$ is even.
Changing a generator of $\cyc{m}$ if necessary we may achieve that
$\wt(x_1)\equiv 1\mod m$. The rest follows from normalizedness.
Taking into account Lemma \ref{even}, from now we may assume that
$\ord(x_1^m)>3\m$, and by symmetry,
$\ord(x_2^m)>3\m$. We will derive a contradiction with $a_3>\m$.

\subsection{}
\label{oc}
By Corollaries \ref{grw} and \ref{formula}, there exist  different
simple
 invariant monomials $\psi_1$, $\psi_2$, $\psi_3$
such that
$$
3\bar{m}\ge i_P\m\ge \mt{ord}(\psi_1)-a_1+\mt{ord}(\psi_2)-a_2
 +\mt{ord}(\psi_3)-a_3-\bar{m}w_P.
$$
Remind that the last inequality can be rewritten as
$$
3\m\ge i_P\bar{m}\ge \mt{ord}(\nu_1)+\mt{ord}(\nu_2)
 +\mt{ord}(\nu_3)-\bar{m}w_P,
$$
where $\nu_i x_i=\psi_i$, $i=1,2,3$.
In particular $\mt{ord}(\nu_1)+\mt{ord}(\nu_2)\le 3\m$.
Since $\wt(\nu_1)\equiv\wt(x_2)\mod m$, $\wt(\nu_2)\equiv\wt(x_1)\mod m$,
by \ref{normal}, we have $a_1+a_2\le 3\m$.
Consider the following cases.

\subsection{Case $a_1+a_2=3\m$.}
Then we have $\ord(\nu_1)+\ord(\nu_2)=3\m$.
Therefore $i_P=3$, $w_P<1$ (by the inequality in \ref{grw}) and
$\ord(\nu_3)=\m w_P<\m$ by \ref{oc}. Since $a_3>\m$,
the monomial $\nu_3$ does not depend on $x_3$.
Up to permutation of $x_1$, $x_2$ we may assume that $\nu_3=x_1^\alpha$
for some $\alpha\in\NN$ such that $\alpha\wt(x_1)+
\wt(x_3)\equiv 0\mod m$ and
$\alpha a_1<\m$. Thus $\psi_3=x_1^\alpha x_3$ and
$a_1<\m$, $a_2>2\m$.
\par
 Further from normalizedness and from $\wt\nu_1\equiv\wt x_2\mod m$,
we have $\ord(\nu_1)=a_2>2\m$, so $\ord(\nu_2)=a_1<\m$.
It gives us $\nu_2=x_1$, $\psi_2=x_1x_2$ and
$\ord(\psi_1)=3\m$. By the assumption in \ref{assumptions1},
$\psi_1$ is not a power of $x_1$. If $\nu_1$ depends on $x_2$,
then $\psi_1=x_1x_2=\psi_2$, a contradiction.
Hence $\psi_1=x_1^\beta x_3^\gamma$, where $\beta,\gamma>0$,
$a_1 \beta + a_3 \gamma=3\m$ and
$\wt(x_1)\beta +\wt(x_3)\gamma\equiv 0\mod m$.
Since $a_3>\m$, $\gamma\le 2$.
On the other hand if $\gamma=1$, then $x_1^{\beta -\alpha}=\psi_1/\psi_3$
is an invariant
of order $\le 2\m$,
a contradiction with assumptions in \ref{assumptions1}.
Thus $\gamma=2$. It gives us $a_1\beta +2a_3=3\m$ and $a_1\alpha+a_3=2\m$.
Then $(2\alpha-\beta )a_1=\m$, so $a_1=1$ and $2\alpha-\beta =\m$,
because $(a_1,\m)=1$. On the other hand
 $\beta \wt(x_1)+2\wt(x_3)\equiv 0\mod m$ and
$\alpha\wt(x_1)+\wt(x_3)\equiv 0\mod m$, so
$\m\wt(x_1)=(2\alpha-\beta)\wt(x_1)\equiv 0\mod m$, the contradiction with
$(m,\wt(x_1))=1$.
\subsection{}
Thus  $a_1+a_2\le 2\m$.
Then there are no invariant  monomials in $x_1, x_2, x_3$
of order $<\ord(x_1x_2)$.
Up to permutation of $\psi_1$ and $\psi_2$ we can take $\psi_1=x_1x_2$.
Since $\psi_2$ is simple, it depends only on $x_2, x_3$.

\subsection{Case $a_1+a_2=2\m$. }
Simple invariant monomials of order $\le 2\m$ may be only  of the form
$x_1x_2$ or $x_i^\alpha x_3$, $i=1,2$\quad
(if $\alpha\wt(x_i)+\wt(x_3)\equiv 0\mod m$
and $\alpha a_i+a_3=2\m$).
On the other hand by the inequality in \ref{oc},
$$
\ord(\psi_2)\le 3\m+a_1+a_2-\ord(\psi_1)=3\m.
$$
There are two subcases.

\subsubsection{Subcase $\ord(\psi_2)=2\m$.}
\label{jjjj}
Then $\psi_2=x_2^\alpha x_3$, because $\psi_2$ depends on $x_2$.
Therefore $a_3<2\m$, $a_2\alpha<\m$ and $a_1>\m$.
Further, if we take in  \ref{computation-w}
$\psi=x_2^\alpha$, we obtain $\m w_P=\ord(x_2^\alpha)=a_2\alpha<\m$.
It follows from  inequality \ref{oc} that
$\ord(\psi_3)\le \m w_P+ a_3+\m\le 3\m$, hence $\ord(\psi_3)=3\m$.
Using $\partial\psi_3/\partial x_3\ne 0$ we get
 $\psi_3\in\{ x_1x_3,\ x_2^\gamma x_3,\ x_2^\gamma x_3^2\}$.
But if $\psi_3=x_1x_3$, then $\wt(x_3)\equiv \wt(x_2)\mod m$,
and by normalizedness $a_3=a_2<\m$, a contradiction.
\par
If $\psi_3=x_2^\gamma x_3$, then
$\ord(x_2^{\gamma-\alpha})=\ord(\psi_3/\psi_2)=\m$.
We derive a contradiction with assumptions in \ref{assumptions1}.
\par
Thus $\psi_3=x_2^\gamma x_3^2$, where
$\gamma \wt(x_2)+2\wt(x_3)\equiv 0\mod m$
and $\gamma a_2+2a_3=3\m$. Consider the monomial
$\nu=x_2^{\alpha-\gamma}=x_2\psi_2/\psi_3$. We have $\ord(\nu)=a_3-\m<a_3$.
On the other hand
$\wt(\nu)\equiv\wt(x_3)\mod m$.
Hence $\wt(\nu)\equiv \wt(x_3)\mod m$ and $\ord(\nu)<\ord(x_3)$,
a contradiction with normalizedness.

\subsubsection{Subcase  $\mt{ord}(\psi_2)=3\m$. }
By the inequality in \ref{oc}, $i_P=3$, $\ord(\psi_3)=a_3+\m w_P$ and by
Corollary \ref{grw}, $w_P<1$. Therefore $\ord(\psi_3)<a_3+\m$.
Thus one has $\psi_2\in \{x_2^\beta x_3^2,\quad x_2^\beta x_3\}$.
But if $\psi_2=x_2^\beta x_3^2$, then $a_2\beta+2a_3=3\m$.
Hence $a_3<2\m$, $a_2<\m$, $a_1>\m$ and  $\ord(\psi_3)=2\m$.
So $\ord(\psi_3)=2\m$ and $\psi_3$ depends on $x_2, x_3$.
Then we can permute $\psi_2$, $\psi_3$ and get the subcase above.
 \par
Therefore $\psi_2=x_2^\beta x_3$. Then $a_2\beta+a_3=3\m$, $a_3<3\m$.
We claim that $a_2<\m$, $a_1>\m$.
Indeed if otherwise, we have $\beta=1$, $\wt(x_3)\equiv\wt(x_1)$ and
by normalizedness $a_3=a_1<\m$, a contradiction.
Further $\ord(\psi_3)\le 3\m$ gives us
$$
\psi_3\in\{x_2^\gamma x_3, \quad x_2^\gamma x_3^2, \quad x_1x_3\}.
$$
 If $\psi_3=x_2^\gamma x_3$, then
the monomial $\psi_2/\psi_3=x_2^{\beta-\gamma}$ is invariant and has
order $\le \m$, a contradiction with assumptions in \ref{assumptions1}.
But if $\psi_3=x_2^\gamma x_3^2$, then we can permute
$\psi_2$, $\psi_3$ and derive a contradiction as above.
Therefore $\psi_3$ depends on $x_1, x_3$, so
$\psi_3=x_1x_3$. Again from normalizedness we have $a_3=a_2<\m$,
a contradiction.

\subsection{Case $a_1+a_2=\m$.} Recall that $\psi_1=x_1x_2$.
As above by the inequality in \ref{oc}
$$
\mt{ord}(\psi_2)\le 3\m,
$$
and if the equality holds, then $i_P=3$, $w_P<1$ and
$$
\ord(\psi_3)=a_3+\m w_P<a_3+\m.
$$
Since $\psi_1=x_1x_2$ and $\psi_2$ are different simple invariant monomials,
we have that $\psi_2$ depends on $x_3$.
Then obviously $\ord(\psi_2)\ge 2\m$.
Consider subcases.

\subsubsection{Subcase $\ord(\psi_2)=2\m$.}
Then it is easy to see  $\psi_2=x_2^\beta x_3$.
It implies that
$$
a_2\beta+a_3=2\m,\qquad \wt(x_1)\beta\equiv \wt(x_3)\mod m.
$$
Moreover we may assume that $\beta\ge 2$ (otherwise we derive 
a contradiction
as in \ref{jjjj}).
 If we take in  \ref{computation-w}
$\psi=\psi_2/x_3=x_2^\beta$, we get
$w_P<1$. Thus $\ord(\psi_3)\le \m+\m w_P+a_3<3\m+\m w_P$ and
$\ord(\psi_3)\le 3\m$.
We claim that $\ord(\psi_3)=3\m$.
Indeed, if otherwise, we have $\psi_3=x_1^\gamma x_3$,
where $a_1\gamma +a_3=2\m$.
It gives us $a_1(\gamma -1)+a_2(\beta-1)+2a_3=3\m$.
Since $\beta\ge 2$, $\gamma =1$. But then $x_1x_3$ is an invariant
and $a_3=a_2<\m$, a contradiction.
\par
Thus $\ord(\psi_3)=3\m$.
Taking into account that $\psi_3$ and $\psi_3/\psi_2$ are
 not powers of $x_i$, $i=1,\dots ,4$, we obtain
$$
\psi_3\in\{ x_1^\gamma x_3,\quad x_1^\gamma x_3^2, 
\quad x_2^\gamma x_3^2\}.
$$
If $\psi_3= x_1^\gamma x_3$, then
 $a_1\gamma +a_3=3\m$, so
$a_3=3\m-a_1\gamma $ and $a_1(\gamma +\beta)=\m (\beta+1)$.
Since $(a_1,\m)=1$, $\gamma +\beta=\m k$, where $k$ is a
positive integer.
Then $\beta=a_1\m k-1\ge\m-1$.
The equality $a_2\beta+a_3=2\m$ is possible only if
$a_1=a_2=1$, $\m=2$, but then $\beta=1$, a contradiction.
\par
Assume that $\psi_3= x_1^\gamma x_3^2$.
Then $a_1\gamma +2a_3=3\m$, $a_1(\gamma -1)+a_2(\beta-1)+3a_3=4\m$.
Since $\beta\ge 2$, $\gamma =1$.
Whence
$$
3\m=a_1+2a_3=a_1+2(2\m-a_2\beta)=4\m+a_1-2a_2\beta=
(4-2\beta)\m+(1+2\beta)a_1.
$$
 Thus $a_1(2\beta+1)=\m(2\beta-1)$.
Since $(a_1,\m)=1$, it implies $\m=2\beta+1$, $a_1=2\beta-1$.
On the other hand $\wt(x_1)+2\wt(x_3)\equiv 0\mod m$ and
$\wt(x_1)\beta\equiv \wt(x_3)\mod m$ imply $\m=1+2\beta\equiv 0\mod m$.
So $\wt(x_1)(1+2\beta)\equiv 0\mod m$, a contradiction.
\par
Finally assume that $\psi_3= x_2^\gamma x_3^2$.
Then $a_2\gamma +2a_3=3\m$ and $a_2\beta+a_3=2\m$.
It gives us $\m=a_2(2\beta-\gamma )$, so $a_2=1$, $\m=2\beta-\gamma $.
On the other hand, $\wt(x_1)\gamma \equiv 2\wt(x_2)\mod m$,
because $\psi_3$ is an
invariant. Whence
$2\wt(x_1)\beta\equiv \wt(x_1)\gamma \mod m$
and $\m=2\beta-\gamma \equiv 0\mod m$,
a contradiction.

\subsubsection{Subcase $\ord(\psi_2)=3\m$.}
Since $\psi_2$ depends on $x_2$ and $x_3$,
$a_3<3\m$. We claim that $a_3<2\m$. Indeed, if otherwise, we have
 $\psi_2= x_2^\beta x_3$
and $\ord(\psi_3)=3\m$, hence $\psi_3= x_1^\alpha x_3$.
Thus $a_1\alpha +a_3=3\m$ and $a_2\beta+a_3=3\m$ give us
$(\alpha-1)a_1+(\beta-1)a_2+2a_3=5\m$.
So $\alpha=1$ or $\beta=1$. But then $x_1x_3$ or $x_2x_3$ is an
invariant. By normalizedness $a_3<\m$, a
  contradiction.
\par
Therefore
$a_3<2\m$, hence $\ord(\psi_3)=2\m$,
$$
\psi_3\in\{ x_1^\gamma x_3,\quad x_2^\gamma x_3
\},
\qquad
\psi_2\in\{x_2^\beta x_3,\quad x_2^\beta x_3^2\}.
$$
 But then
we can permute $\{\psi_2,\ \psi_3\}$ and $\{x_1,\ x_2\}$ and get the
 previous subcase.
This proves Lemma \ref{matrix}.
\qq

\subsection{Proof of Theorem \ref{main}.}
\label{lll}
First we treat  \ref{matrix2} (cf. \cite[(4.4)]{Mori-flip}).
We claim that  in case (i) of \ref{matrix2} $(X,P)$ is a 
cyclic quotient singularity.
 Let $C\3=C\3(1)\cup C\3(2)$ be irreducible decomposition of
$C\3=(\pi^{-1}(C))_{\mt{red}}$.
The component $C\3(1)$ is the image of
$t\longrightarrow (t, t^{\m-1}, t^{\m+1}, t^{\m})$. The curve 
$C\3(1)$ is
a complete intersection:
$$
C\3(1)=\{ x_1x_2-x_4=x_3-x_1^{\m+1}=x_2-x_1^{\m-1}=0 \}
$$
and so is $C\3$ :
$$
C\3=\{ x_1x_2-x_4=x_3-x_1^{\m+1}=x_2^2-x_1^{2\m-2}=0 \}.
$$
Therefore the equation of $X\3$ can be written as
$$
\phi=(x_1x_2-x_4)\varphi_1+
(x_3-x_1^{\m+1})\varphi_2+
(x_2^2-x_1^{2\m-2})\varphi_3=0,
$$
where $\varphi_i=\varphi_i(x_1,x_2,x_3,x_4)$ are semi-invariants.
Since $\wt(x_3)\not\equiv \pm \wt(x_1),\wt(x_2)$ and because
$(X,P)$ is of type $cA/m$, the equation
of $\phi$ must contain either $x_1x_2$ or $x_4$ (see \cite{Mori-term}).
But this is possible only if $\varphi_1(0,0,0,0)\ne 0$.
Thus $\phi$ always contains $x_4$ and $X\3$ is non-singular.
We obtain case (iv) of \ref{main1}.
Case  (ii) of \ref{matrix2} is treated by the similar way.
\subsection{}
\label{iiii}
From now we will suppose that $f$ is not from \ref{main1} (iv) or (v).
Consider the base change from \ref{covering}
$$
\begin{CD}
X\6 @>{g}>>& X\\
@V{f\6}VV @V{f}VV\\
S\6 @>{h}>> S
\end{CD}
$$
If $\m=1$, then $X\6$ is Gorenstein and $f\6$ is a "classical" conic bundle
\cite{Cut} (it means that $S\6$ is non-singular and
 every fiber of $f\6$ is a conic in $\PP^2$).
Since $(X,C)$ is imprimitive,
 the central fiber $C\6$ has exactly two components.
The group $\cyc{2}$ permutes
these  components. In this case it is not difficult to
write down an equation of $X\6$ explicitly (see \cite{Pro2} or \cite{Pro3}).
We obtain case \ref{main1}, (i).
\par
Assume now that $\m>1$.
By Lemma \ref{matrix},
the general member $F\in |-K_X|$ does not contain $P$
and has DuVal point at $P$.
The cover $f_F\colon (F,P)\to (S,0)$ is finite of degree 2.
Thus $(S,0)$ is a quotient of $(F,P)$ by an involution.

\begin{proposition1}
\textup{(} see \cite{Cat}\textup{)}
 \label{catanese}
Let $(F,P)$ be a germ of DuVal singularity and
$\tau \colon (F,P)\to (F,P)$ be an  (analytic) involution. Then
up to analytic isomorphism there are only the following
possibilities.
\label{catanese2}
\par\bigskip
$$
\begin{array}{|c|c|c||c|c|c|}
\hline
\mt{no.}&(F,P)&(F,P)/\tau&\mt{no.}&(F,P)&(F,P)/\tau\\
\hline
&\qquad&\qquad\qquad\qquad&&\qquad&\qquad\\
1) &\hbox{{\rm any}}&\hbox{{\rm smooth}}         &
6) &A_{2k+1}        &\CC^2/\cyc{4k+4}(1,2k+1)    \\
2) &A_{2k+1}        &D_{k+3}                     &
7) &E_6             &A_2                         \\
3) &E_6             &E_7                         &
8) &A_{2k}          &\CC^2/\cyc{2k+1}(1,2k-1)    \\
4) &D_k             &D_{2k-2}                    &
9) &D_k             &A_1                         \\
5) &A_k             &A_{2k+1}                    &
10)&A_{2k+1}        &A_k                         \\
\hline
\end{array}
$$
\end{proposition1}
\par\bigskip
\begin{remark3}
In cases 6) and 8) the dual graphs of corresponding minimal
resolutions
of $(F,P)/\tau$ are
$$
\begin{array}{c}
{\scriptstyle -4}\\
\circ\\
\end{array}
\qquad
(\hbox{{\rm for}}
\quad k=0),
\leqno 6)
$$
$$
\underbrace{
\begin{array}{ccccccc}
{\scriptstyle -3}&&{\scriptstyle -2}&
&{\scriptstyle -2}&&{\scriptstyle -3}\\
\circ&\pal&\circ&\cdots&\circ&\pal&\circ\\
\end{array}}_{k+1}
\quad
(\hbox{{\rm for}}
\quad k>0)
\leqno 6)
$$
$$
\underbrace{
\begin{array}{ccccccccc}
{\scriptstyle -2}&&{\scriptstyle -2}&
&{\scriptstyle -2}&&{\scriptstyle -2}&
&{\scriptstyle -3}\\
\circ&\pal&\circ&\cdots&\circ&\pal&\circ&\pal&\circ\\
\end{array}}_{k}
\leqno 8)
$$
As noticed by Shokurov, these singularities are exactly all log-terminal
singularities with discrepancies $\ge -1/2$.
\end{remark3}
\subsubsection{Proof of Theorem \ref{main} (continue).}
Since $(S,0)=(F,P)/\tau$ is a cyclic quotient singularity  by
Proposition \ref{ooo},
cases 1)-4) are impossible.
Further by assumptions in \ref{iiii} we have $\m>1$, hence $m>d$.
Remind that in our situation $d$ is the order of local fundamental
group $\pi_1(S-\{0\})$ (see \ref{odnatochka}).
In particular, $d=n+1$ if $(S,0)$ is a point of type $A_n$.
Thus in case 5) we have $m>d=2k+2\ge 4$.
But $(F,P)$ is a point of type $A_k$ and has the topological index $k+1$.
This contradicts to \ref{index-surf-threef}. The same arguments show that
the cases 6),7),8) are also impossible (if $\m>1$).
In case 9) $d=2$. From the list of \ref{term-cl} we see that
$m>2$ only if $(X,P)$ is of type $cAx/4$ and then $\m=2$, $m=4$.
This is case (ii) of  \ref{main1} by Lemma \ref{except}.
\par
Finally consider case 10). Then $d=k+1$.
In the classification \ref{term-cl} $(X,P)$ is of  type
$cA/n$. Point $(F,P)$ has topological index $2k+2$ and $d=k+1$, so by
Lemma \ref{index-surf-threef}
$(k+1)\m=d\m=m\le 2(k+1)$, whence  $\m=2$.
Moreover $(X,P)$ is a cyclic quotient singularity in this case
(because an equality in \ref{index-surf-threef} holds).
Lemma \ref{warwick} gives us then $d=2$ or $4$.
\par
Suppose that $d=2$ (and $\m=2$).
From Lemma \ref{matrix}, $a_3=1$. By changing a character $\chi$
if necessary we get $\wt(x_3)\equiv 1\mod 4$. Up to permutation $\{1,\ 2\}$
we also have $\wt(x_1)\equiv 1\mod 4$ and $\wt(x_2)\equiv 3\mod 4$.
Then by normalizedness, $a_1=1$ and $a_2=1$ or $3$.
Obviously $\wt(x_4)\equiv 0\mod 4$ and $a_4=\m=2$.
Thus a component $C\3(1)$ of $C\3$ is the image of
$t\longrightarrow (t,t^{a_2},t,t^2)$,
where $a_2=1$ or $3$. But if $a_2=3$, then $C\3$ can be given by
$x_2-x_1^3=x_3-x_1=x_4^2-x_1^4=0$. Whence
an equation of $X\3$ in $\CC^4$ in a normalized coordinate system
 is
$\phi=(x_2-x_1^3)\varphi_1+(x_3-x_1)\varphi_2+(x_4^2-x_1^4)\varphi_3=0$.
Since $(X,P)$ is a cyclic quotient, $\phi$ must contain the term $x_4$,
a contradiction. Therefore $a_2=1$ and  $C\3$ is given by
$x_2^2-x_1^2=x_3-x_1=x_4-x_1x_2=0$. By changing coordinates
we obtain \ref{main1}(iii).
\par
If $d=4$, then similarly one can see
$a_3=1$ and $\wt(x_3)\equiv 1\mod 8$. Again an equation of $X\3\subset\CC^4$
in a normalized coordinate system
contains the term $x_4$, so it also contains a monomial $\phi_0\ne x_4$
such that $\wt(\phi_0)\equiv 0$ and $\ord(\phi_0)=\ord(x_4)=2$.
Then $\phi_0=x_1x_2$ and $a_1=a_2=1$. Denote $a:=\wt(x_1)$.
It is easy to see that we can take $a\in \{1,\ 3\}$.
If $a=1$, then $C\3=\{ x_2^4-x_1^4=x_3-x_1=x_4-x_1x_2=0\}$.
If $a=3$, then one can check that
$C\3=\{ x_2^2-x_3^2=x_2x_3-x_1^2=x_4-x_1x_2=0\}$.
After changing coordinates  we obtain  \ref{main2}.
\par
To complete the proof we have to show that $X$ can contain at most one
Gorenstein singular point if $(S,0)$ is of type $A_1$ and
contains no Gorenstein singular point if $(S,0)$ is of type $A_3$.
Indeed in diagram \ref{covering}
$f\6\colon (X\6,C\6)\to (S\6,0)$ is a Mori conic bundle with reducible
central fiber $C\6$. Therefore by \ref{vanishing}, $\rho(X\6/S\6)>1$ and
each component $C\6(i)$ of $C\6$ generates an extremal ray on $X\6$
over $(S\6,0)$. Thus each $(X\6,C\6(i))$ is an (not necessary flipping)
extremal neighborhood in the sense of \cite{Mori-flip}.
By \cite[Theorem (6.2)]{Mori-flip}, we see that  $(X\6,C\6(i)$
contains at most one Gorenstein singular point and so does $(X,C)$
because $X\6\to X$ is \'etale outside $P$.
\par
In case \ref{main2} $X\6$ has index 2, whence
$(-K_{X\6}\cdot C\6(i))\in \frac{1}{2}\ZZ$.
On the other hand  as in \ref{xaxa} $f\6\colon (X\6,C\6)\to (S\6,0)$ is flat.
Therefore for a general fiber $L$ of $f\6$ we have $L\equiv r\sum C\6(i)$,
where $r\in\NN$. Thus $2=(-K_{X\6}\cdot L)=4r(-K_{X\6}\cdot C\6(i))$, so
$1/2r\in\frac{1}{2}\ZZ$. It gives us $r=1$, $L\equiv \sum C\6(i)$.
Whence the scheme-theoretical central fiber of $f\6$ over $0$ is
generically reduced. Then $C\6$ is a complete intersection of Cartier
 divisors inside
$X\6$. Since $C\6$ has a unique singular point $P\6$, $X\6$ has exactly one
singular point (at $P\6$) along $C\6$ and so has $X$ (at $P$).
In case (i) of \ref{main1} one  can use the same arguments to show
that the  singular point $P$ is unique.
This proves Theorem \ref{main}.
\qq
\par
 Finally we propose examples of Mori conic bundles over a
non-singular base surface (case \ref{main0} of Theorem).

\begin{example}
Let $Y\subset\PP^3_{x,y,z,t}\times\CC^2_{u,v}$
is given by the equations

$$
\left\{
\begin{array}{l}
xy-z^2=ut^2\\
x^2=uy^2+v(z^2+t^2)
\end{array}
\right.
$$
It is easy to check that $Y$ is non-singular.
The projection $\PP^3\times\CC^2\to\CC^2$ gives us an elliptic fibration
$g\colon Y\to\CC^2$. A general fiber of $g$ is an 
intersection of two quadrics in
$\PP^3$ and the central fiber $g^{-1}(0)$ is multiple $\PP^1$.
  Let $\cyc{2}$ acts on $Y$ by
$$
x\to -x,\qquad
y\to -y,\qquad
z\to -z,\qquad
t\to t,\qquad
u\to u,\qquad
v\to v.
$$
The locus of fixed points consist of the divisor $B:=\{ t=0\}\cap Y$
and the isolated point $(x,y,z,t;u,v)=(0,0,0,1;0,0)$.
Consider the quotient mrphism $h\colon Y\to X:=Y/\cyc{2}$.
By the Hurwitz formula $0=K_Y=h^*(K_X)+D$, where $h^*(-K_X)=D$.
Whence $-K_X$ is ample over $\CC^2$.
Therefore the quotient $f\colon X\to \CC^2$ is a Mori conic bundle with
irreducible central fiber that contains
 only one singular point of type $\frac{1}{2}(1,1,1)$.
\end{example}
 The following example shows that the total space $X$
of a Mori conic bundle can contain two Gorenstein terminal
points (cf. \cite[0.4.13.1]{Mori-flip}).

\begin{example}
Let $Y\subset\PP^3_{x,y,z,t}\times\CC^2_{u,v}$
is given by the equations
$$
\left\{
\begin{array}{l}
x^2=uz^2+vt^2\\
y^2=ut^2+vz^2
\end{array}
\right.
$$
As above the quotient $X:=Y/\cyc{2}(1,1,1,0;0,0)\to \CC^2$
is a Mori conic bundle with
irreducible central fiber.
Singular points  of $X$ are one point of
type $\frac{1}{2}(1,1,1)$ and two ordinary double points.
\end{example}

\bibliographystyle{amsalpha}

\end{document}